\documentclass[11pt,a4paper]{article}
\usepackage{latexsym}
\usepackage{epsfig}
\usepackage{here}
\usepackage{amsmath}
\usepackage{amssymb}
\usepackage{amsfonts}
\usepackage{amsbsy}
\usepackage{graphicx}
\usepackage{graphics}
\usepackage{euscript}
\usepackage{epsfig}
\usepackage{varioref}
\usepackage{enumerate}
\usepackage{verbatim}
\usepackage{here}
\usepackage{makeidx}
\usepackage{color}

\makeindex

\def\R{\mathbb{R}}

\def\P{\mathbb{P}}

\def\N{\mathbb{N}}
\def\u{{\bf u}}
\def\v{{\bf v}}
\def\w{{\bf w}}

\def\x{{\bf x}}
\def\f{{\bf f}}
\def\a{{\bf a}}

\newtheorem{theorem}{Theorem}[section]

\newtheorem{remark}{Remark}[section]
\newtheorem{corollary}{Corollary}[section]
\newtheorem{lemma}{Lemma}[section]
\numberwithin{equation}{section} 

\makeatletter
\renewcommand\section{\@startsection {section}{1}{\z@}%
 				   {-2.5ex \@plus -1ex \@minus -.2ex}%
                                   {1.3ex \@plus.2ex}%
                                   {\normalfont\large\bfseries}}
\renewcommand\subsection{\@startsection {subsection}{1}{\z@}%
 				   {-2.ex \@plus -1ex \@minus -.2ex}%
                                   {0.5ex \@plus.2ex}%
                                   {\normalfont\normalsize\bfseries}}
\makeatother

\begin{document}

\title{$L^2$-stability of explicit schemes for incompressible Euler equations}

\author{
Erwan Deriaz \thanks{Institute of Mathematics, Polish Academy of Sciences. ul. Sniadeckich 8, 00-956 Warszawa,
Poland, to whom correspondance should be addressed (E.Deriaz@impan.gov.pl)}
}


\bigskip
\bigskip

\date{\today}

\maketitle

\begin{abstract}

We present an original study on the numerical stabiliy of explicit schemes solving the incompressible
Euler equations on an open domain with slipping boundary conditions.
Relying on the skewness property of the non-linear term, we demonstrate that some explicit schemes
are numerically stable for small perturbations under the condition $\delta t\leq C \delta x^{2r/(2r-1)}$
where $r$ is an integer, $\delta t$ the time step and $\delta x$ the space step.

\end{abstract}

\section{Introduction}
\label{intro}

In order to achieve stability of a numerical scheme solving the incompressible Euler equations in a divergence-free
dicretisation frame, we seek for a criterion that link the time and the space steps.
During unsteady incompressible fluid simulations with the help of divergence-free wavelets \cite{DP07}, we observed
the CFL-like condition: $\delta t\leq C\delta x^{4/3}$ for a centered scheme order two (\ref{SC2}), presented in \cite{KT97,DP07}.
Either bibliographical data don't give a satisfying explaination to this phenomena  \cite{Tem84,MT98}, either it relies
on a Von Neumann stability analysis for the advection equation where velocity is a constant vector \cite{Wes01}, which is
not our case.

\ 

This work takes advantage of ideas presented in R.~Temam \cite{Tem84,MT98} for its context: incompressible Euler equations
with the use of the skewness property (lemma \ref{lemme}).
But, at the difference with references \cite{Tem84,MT98}, we are not in the context of finite elements,
but assume a divergence-free space of discretisation.
Meanwhile, we make the connection with the results of Von Neumann stability for the convection equation \cite{Wes01}.

In the following Note, we manage to derive the propagation in $L^2$ norm of a small
perturbation $\varepsilon_n$ of a regular solution $\u_n$ of the numerical explicit scheme \ref{glxplschm} to
solve the incompressible Euler equations (\ref{Eulervp}) in a relatively accurate way.

\ 

In this case, we establish a stability criterion of the form $\delta t\leq C\delta x^{2r/(2r-1)}$ with $r$ an integer
and $C$ an explicitly computable constant, for various explicit schemes.
Contrarily to Von Neumann stability \cite{Wes01} which assumes infinite or periodic domain and constant
advection velocity, we establish this stability condition solely under regularity assumptions on the
velocity.

\ 

The Euler equations modelise incompressible fluid flows with no viscous term:
\begin{eqnarray}
\label{Eulervp}
\frac{\partial \mathbf{u}}{\partial t}+({\bf u}\cdot\nabla) {\bf u}-\nabla p=0,\qquad 
{\rm div}\,\u=0
\end{eqnarray}
The use of the Leray projector $\mathbb{P}$ which is the $L^2$-orthogonal projector on the divergence-free
space, allows us to remove the presure term:
\begin{eqnarray}
\label{EulervPdiv}
\frac{\partial \mathbf{u}}{\partial t}+\mathbb{P}\left[({\bf u}\cdot\nabla) {\bf u}\right]=0
\end{eqnarray}
We start from a discretisation $\u_n$ of the solution $\u$ in time and in space. Then we consider a
perturbed solution $\u_n+\varepsilon_n$. Assuming the regularity of $\u$ and the consistency of
the scheme, we are interested in the evolution of the perturbation $\varepsilon_n$ for different
explicit schemes in time.

\section{Discretisation in time and in space}
\label{discret}

In order to solve equation (\ref{EulervPdiv}) numerically, we discretise the solution
$\u_n=\u(n\delta t)$ in a divergence-free space $V_{{\rm div}\,0}(\delta x)$.
Such spaces appear in spectral codes using the Fourier transform \cite{CHQZ88},
or are produced in a stable way thanks to divergence-free wavelets \cite{U00,DP07}.
The parameter $\delta x$ stands for the smallest space step in the discretisation space
$V_{{\rm div}\,0}(\delta x)$.
Hence, every function $\u_n$ in $V_{{\rm div}\,0}(\delta x)$ satisfies:
$$
\|\partial_i\u_n\|_{L^p}\leq C(p)\delta x^{-1}\|\u_n\|_{L^p}
$$
where $1\leq p\leq +\infty$, and $\partial_i \u_n$ denotes a partial space derivative of $\u_n$.
The constant $C(p)$ can often be taken equal to $1$, and we will do so in the following.

Concerning the discretisation in time, we consider different explicit schemes proceeding
in several steps, like Runge-Kutta schemes:
\begin{equation}
\label{glxplschm}
\u_{n(0)}=\u_n,\quad
\u_{n(\ell)}=\sum_{i=0}^{\ell-1}a_{\ell i}\,\u_{n(i)}-\sum_{i=0}^{\ell-1}b_{\ell i}\,\delta t \,\widetilde{\P}\left[(\u_{n(i)}\cdot\nabla)\u_{n(i)}\right]
~~{\rm for}~1\leq \ell \leq k,\quad
\u_{n+1}=\u_{n(k)}
\end{equation}
where $\widetilde{\P}$ stands for the orthogonal projector on the discretisation space $V_{{\rm div}\,0}(\delta x)$.
One can notice that $\widetilde{\P}\circ\P=\P\circ\widetilde{\P}=\widetilde{\P}$.

\section{$L^2$-stability condition for a small perturbation}
\label{scfl}

Actually, the small error $\varepsilon_n$ that we introduce corresponds to oscillations at the smallest scale
in space $V_{{\rm div}\,0}(\delta x)$. This stability error propagates and may increase at each time step.
In what follows, we demonstrate that under some precise CFL-like conditions, the $L^2$ norm of
this small error $\varepsilon_n$ is amplified such that:
\begin{equation}
\label{stabCdt}
\|\varepsilon_{n+1}\|_{L^2}\leq (1+C\delta t) \|\varepsilon_{n}\|_{L^2}
\end{equation}
where $C$ is a constant that neither depends on $\delta x$ nor on $\delta t$.

Thus, after a time elapse $T$, the error increases at most exponentially as a function of the time:
\begin{eqnarray}
\|\varepsilon_{t_0+T}\|_{L^2}\leq
\left(1+C\delta t\right)^{T/\delta t} \|\varepsilon_{t_0}\|_{L^2}
\leq e^{CT} \|\varepsilon_{t_0}\|_{L^2}
\end{eqnarray}

\ 

For the following stability study, we will need the skewness property of the transport term.
This property is utilized for the stability of the incompressible Navier-Stokes equations
in \cite{MT98}.

\begin{lemma}
\label{lemme}
Let $\u,\v,\w\in H^1(\Omega)^d$, $H^1(\Omega)$ denoting the Sobolev space on the open set $\Omega\subset \R^d$,
be such that $(\u\cdot\nabla)\v,(\u\cdot \nabla)\w \in L^2$.
If $\u\in {\textbf{H}}_{{\rm div},0}(\Omega)=\{\f\in (L^2(\Omega)))^d,~{\rm div}~\f=0\}$, then
$$
<\v,(\u\cdot \nabla)\w>_{L^2(\Omega)}=-<(\u\cdot \nabla)\v,\w>_{L^2(\Omega)}
$$
\end{lemma}

\begin{corollary}
\label{corollary}
With the same assumptions as in lemma \ref{lemme},
$$
<\v,(\u\cdot \nabla)\v>_{L^2(\Omega)}=\int_{\x\in\Omega}\v\cdot(\u\cdot \nabla)\v\, d\x=0
$$
\end{corollary}

\ 

In the scheme (\ref{glxplschm}), we denote by $\varepsilon_{n(\ell)}$ the stability error at level $\ell$.
Then, under the condition $\delta t=o(\delta x)$ and for $\varepsilon_{n}$ small enough,
most of the terms appearing in the expression of $\varepsilon_{n(\ell)}$ become
negligeable compared with:
\begin{itemize}
\item the terms $\delta t^i\,\widetilde{F}^i(\varepsilon_{n})$ where
$\widetilde{F}(\varepsilon_n)=\widetilde{\P}[(\u_n\cdot\nabla)\varepsilon_n]$ and
$\widetilde{F}^i=\widetilde{F}\circ \widetilde{F}\circ \dots \circ \widetilde{F}$, $i$ times.
\item the term $\delta t\, \widetilde{\P}[(\varepsilon_{n}\cdot\nabla)\u_n]$,
\end{itemize}
Hence, for the scheme (\ref{glxplschm}), we find the following expression for $\varepsilon_{n(\ell)}$:
\begin{equation}
\label{dvpmt_Fi}
\varepsilon_{n(\ell)}=\sum_{i=0}^{\ell} \alpha_{\ell i}\delta t^i\widetilde{F}^i(\varepsilon_n)+\beta_\ell\delta t\, \widetilde{\P}[(\varepsilon_{n}\cdot\nabla)\u_n]+o()
\end{equation}
where function $o()$ gathers the negligeable terms.
$$
\alpha_{\ell i}=\sum_{j=i}^{\ell-1} a_{\ell j}\alpha_{j i}+\sum_{j=i-1}^{\ell-1} b_{\ell j}\alpha_{j i-1}\quad
{\rm and}\quad \beta_{\ell}=\sum_{j=1}^{\ell-1} a_{\ell j} \beta_{j} + \sum_{j=0}^{\ell-1} b_{\ell j}
$$
As a result, we find the following expression for $\varepsilon_{n+1}$:
\begin{equation}
\label{e_n+1_Fi}
\varepsilon_{n+1}=\sum_{i=0}^{k} \alpha_{i}\delta t^i\widetilde{F}^i(\varepsilon_n)-\delta t\, \widetilde{\P}[(\varepsilon_{n}\cdot\nabla)\u_n]+o()
\end{equation}
Starting from this expression and using the fact that according to lemma \ref{lemme},
\begin{equation}
\label{FF}
<\widetilde{F}^i(\varepsilon_n),\widetilde{F}^j(\varepsilon_n)>_{L^2(\Omega)}=\left\{
\begin{array}{ll} 0 &~~~~{\rm if}~~~~i+j=2\ell+1~~~~{\rm for}~~~~\ell\in\N \\
(-1)^{\ell-i}\|\widetilde{F}^\ell(\varepsilon_n)\|_{L^2}^2 &~~~~{\rm if}~~~~i+j=2\ell~~~~{\rm for}~~~~\ell\in\N
\end{array}\right.
\end{equation}
we compute the $L^2$ norm of $\varepsilon_{n+1}$ as a function of the $L^2$ norm of $\varepsilon_{n}$:
\begin{equation}
\label{sumS}
\|\varepsilon_{n+1}+\delta t\, \widetilde{\P}[(\varepsilon_{n}\cdot\nabla)\u_n]\|_{L^2}^2=\sum_{\ell=0}^{k} S_{\ell}\,\delta t^{2\ell}\|\widetilde{F}^\ell(\varepsilon_n)\|_{L^2}^2+o()
\end{equation}
with
\begin{equation}
\label{Sell}
S_{\ell}=\sum_{j=-{\rm min}(\ell,k-\ell)}^{{\rm min}(\ell,k-\ell)}(-1)^j\alpha_{\ell-j}\alpha_{\ell+j}
\end{equation}
For consistency needs of the numerical scheme, we must have $S_0=1$. If, on an other hand we suppose
$S_1=S_2=\dots=S_{r-1}=0$ and $S_{r}>0$, knowing that in the discretised space $V_{{\rm div}\,0}(\delta x)$,
$$
\|\widetilde{F}^r(\varepsilon_n)\|_{L^2}\leq \|{\bf u}_n\|_{L^\infty}^{r} \frac{\|\varepsilon_n\|_{L^2}}{\delta x^r}
\quad {\rm and}\quad
\|\widetilde{\P}[(\varepsilon_{n}\cdot\nabla)\u_n]\|_{L^2}\leq\|\varepsilon_{n}\|_{L^2}\|\nabla \u_n\|_{L^\infty}
$$
we derive:
\begin{equation}
\label{expS}
\|\varepsilon_{n+1}\|_{L^2}\leq\left( 1+\left(\|\nabla \u_n\|_{L^\infty}+\frac{\delta t^{2r-1}S_r}{2\delta x^{2r}}\|{\bf u}_n\|_{L^\infty}^{2r}+o()\right)\delta t \right) \|\varepsilon_{n}\|_{L^2}
\end{equation}
If, on an other hand we assume the consistency, there exist constants $A_0$ and
$A_1$ such that $\|{\bf u}_n\|_{L^\infty}\leq A_0$ and $\|\nabla \u_n\|_{L^\infty}\leq A_1$
when $\delta x$ and $\delta t$ go to $0$.
Hence, the numerical scheme (\ref{glxplschm}) is stable for small pertubations under the condition:
\begin{equation}
\label{CFL1}
\delta t\leq C \delta x^{\frac{2r}{2r-1}}
\end{equation}
That brings the following theorem out:
\begin{theorem}
\label{theoO_CFL}
An order $2p$ scheme solving the incompressible Euler equations is numerically stable for
small perturbations at worst under the CFL-like condition:
\begin{equation}
\label{CFL2}
\delta t\leq C \delta x^{\frac{2p+2}{2p+1}}
\end{equation}
\end{theorem}
\emph{Proof:} 
For an order $2p$ scheme, we have the following equality, point by point:
$$
\u_{n+1}=\u_n+\delta t\, \partial_t\u_n+\frac{\delta t^2}{2}\partial_{tt}\u_n+\dots
+\frac{\delta t^{2p}}{(2p)!}\partial_t^{2p}\u_n+o(\delta t^{2p})
$$
Considering $\partial_t\u_n=\P[(\u_{n}\cdot\nabla)\u_n]$ and introducing a small perturbation $\varepsilon_{n}$,
leads to:
\begin{equation}
\label{expeps}
\varepsilon_{n+1}=\varepsilon_n+\delta t \, \widetilde{F}(\varepsilon_n)+\frac{\delta t^2}{2} \widetilde{F}\circ F(\varepsilon_n)+\dots
+\frac{\delta t^{2p}}{(2p)!} \widetilde{F}\circ F^{2p-1}(\varepsilon_n)+\delta t \,\widetilde{\P}[(\varepsilon_{n}\cdot\nabla)\u_n]+o()
\end{equation}
with $F(\varepsilon)=\P[(\u_{n}\cdot\nabla)\varepsilon]$ and $o()$ gathering the terms that are negligeable
under the condition $\delta t=o(\delta x)$. Then for $q\in[1,p]$,
$$
S_q=\sum_{p=0}^{2q}(-1)^{(q-p)}\frac{1}{p!}~\frac{1}{(2q-p)!)}=\frac{(-1)^q}{(2q)!}\sum_{p=0}^{2q}{\rm C}_{2q}^p(-1)^p
=0
$$
Which allows us to conclude, as stated at line (\ref{expS}).

\ 

\begin{remark}
\label{VonNeumann}
A Von Neumann stability analysis would have proceeded as follows.
We compute the evolution of the Fourier mode $\varphi(n\delta t)=\varphi_n e^{i\zeta\cdot \x}$ with $\zeta\in \R^d$,
for the advection equation $\partial_t \varphi=-\a\cdot \nabla \varphi$, with $\a=\u$ a constant velocity.
As $\nabla \varphi=i\zeta\cdot\varphi$, for the scheme (\ref{glxplschm}), we find $\varphi_{n+1}=\xi\varphi_{n}$
with, as for computation (\ref{e_n+1_Fi}),
$$
\xi=\sum_{j=0}^{k} \alpha_{j}(-i\a\cdot \zeta)^j\quad {\rm and~then}\quad
|\xi|^2=\sum_{\ell=0}^{k} S_{\ell}\, \delta t^{2\ell} |\a\cdot \zeta|^{2\ell}
$$
The coefficients $S_{\ell}$ have the same expression (\ref{Sell}) as when we used the skewness property.\\
As, on an other hand, we have $|\a\cdot \zeta|\leq \|\a\|/\delta x$ in the discretisation space
$V_{{\rm div}\,0}(\delta x)$, we find the same stability criterion if we want to have $|\xi| \leq 1+C\delta t$.
\end{remark}

\section{Examples}
\label{exemples}

Let $A_0=\sup_{t\in[0,T],~\x\in\Omega}| \u(t,\x)|$ and $A_1=\sup_{t\in[0,T],~\x\in\Omega}| \nabla \u(t,\x)|$.
We propose to apply our stability analysis to some classical schemes.

The simplest example is the Euler explicit scheme, order one in time:
\begin{equation}
\label{EE}
{\bf u}_{n+1}= {\bf u}_{n} - \delta t ~\widetilde{\P}\left[({\bf u}_n\cdot\nabla) {\bf u}_n\right]
\end{equation}
For this, we find:
$$
\| \varepsilon_{n+1} \|_{L^2}\leq\left( 1+ (\frac{A_0^2}{2}\frac{\delta t}{\delta x^2}+A_1)\delta t\right)\|\varepsilon_{n}\|_{L^2},
\qquad {\rm and~the~CFL:}\quad \delta t \leq 2 C \left(\frac{\delta x}{A_0} \right)^2
$$
An improved version of this scheme allows us to construct an order two centered scheme:
\begin{equation}
\label{SC2}
\left\{
\begin{array}{ll}
{\bf u}_{n+1/2}=& {\bf u}_{n} - \frac{\delta t}{2} \widetilde{\P}\left[({\bf u}_n\cdot\nabla) {\bf u}_n\right]\\
\vspace{-0.3cm} \\
{\bf u}_{n+1}=& {\bf u}_{n} - \delta t~ \widetilde{\P}\left[({\bf u}_{n+1/2}\cdot\nabla) {\bf u}_{n+1/2}\right]
\end{array}
\right.
\end{equation}
For this scheme, stability is slightly improved:
$$
\|\varepsilon_{n+1}\|_{L^2}\leq \left( 1+\frac{\delta t^4}{8\delta x^4}A_0^4+\delta t A_1 \right)\|\varepsilon_{n}\|_{L^2}
\qquad \textrm{hence~the~CFL~:}\quad \delta t \leq 2 C^{1/3} \left(\frac{\delta x}{A_0} \right)^{4/3}
$$
For Runge-Kutta scheme of order 4:
{\small
\begin{equation}
\label{RK4}
\left\{
\begin{array}{ll}
{\bf u}_{n(1)}=& {\bf u}_{n} - \frac{\delta t}{2} \widetilde{\P}\left[({\bf u}_n\cdot\nabla) {\bf u}_n\right]\\
\vspace{-0.3cm} \\
{\bf u}_{n(2)}=& {\bf u}_{n} - \frac{\delta t}{2} \widetilde{\P}\left[({\bf u}_{n(1)}\cdot\nabla) {\bf u}_{n(1)}\right]\\
\vspace{-0.3cm} \\
{\bf u}_{n(3)}=& {\bf u}_{n} - \delta t ~\widetilde{\P}\left[({\bf u}_{n(2)}\cdot\nabla) {\bf u}_{n(2)}\right]\\
\vspace{-0.3cm} \\
{\bf u}_{n+1}=& {\bf u}_{n} - \frac{\delta t}{6} \widetilde{\P}\left[({\bf u}_{n}\cdot\nabla) {\bf u}_{n}\right]
- \frac{\delta t}{3} \widetilde{\P}\left[({\bf u}_{n(1)}\cdot\nabla) {\bf u}_{n(1)}\right]
- \frac{\delta t}{3} \widetilde{\P}\left[({\bf u}_{n(2)}\cdot\nabla) {\bf u}_{n(2)}\right]
  - \frac{\delta t}{6} \widetilde{\P}\left[({\bf u}_{n(3)}\cdot\nabla) {\bf u}_{n(3)}\right]
\end{array}
\right.
\end{equation}
}
theorem \ref{theoO_CFL} allows us to predict a CFL-like condition $\delta t\leq C \delta x^{6/5}$ at worst.
Computations show that actually: 
$$
S_1=S_2=0\quad{\rm and}\quad S_3=-\frac{1}{72},\,S_4=\frac{1}{576}
$$
Hence our study doesn't fully apply to this case.

\ 

The order two Adams-Bashford scheme doesn't remain to the definition (\ref{glxplschm}).
It goes as follows:
\begin{equation}
\label{AB}
{\bf u}_{n+1}= {\bf u}_{n} - \frac{3}{2} \delta t ~\widetilde{\P}\left[({\bf u}_{n}\cdot\nabla) {\bf u}_{n}\right]
+\frac{1}{2} \delta t ~\widetilde{\P}\left[({\bf u}_{n-1}\cdot\nabla) {\bf u}_{n-1}\right]
\end{equation}
Nevertheless, it is possible to perform computations similar to those of part \ref{scfl}.
First, as it is order two, according to theorem \ref{theoO_CFL}, it is stable at worst
under the condition $\delta t\leq C\delta x^{4/3}$. Further computations show that:
$$
\|\varepsilon_{n+1}\|_{L^2}\leq \left( 1+\frac{\delta t^4}{4\delta x^4}A_0^4+\delta t A_1 \right)\|\varepsilon_{n}\|_{L^2}
\qquad \textrm{inducing~the~CFL:}\quad \delta t \leq 2^{2/3} C^{1/3} \left(\frac{\delta x}{A_0} \right)^{4/3}
$$

\ 

As a conclusion, many usual schemes for simulating the fluid flows verify a stability
condition of the type $\delta t\leq C \delta x^{2r/(2r-1)}$ with $r$ an integer.

On can remark that if we have relation (\ref{expeps}) at order $m$ with no term of the type
$\widetilde{F}\circ\dots\circ\widetilde{F}(\varepsilon_n)$ in $o()$ (which is, for instance, what Runge-Kutta
schemes satisfy), then the related scheme will have to verify a CFL-like condition of the type $\delta t\leq C \delta x^{(m+1)/m}$
if $m\equiv 1 [4]$ and $\delta t\leq C \delta x^{(m+2)/(m+1)}$ if $m\equiv 2 [4]$.



\section*{Acknowledgements}

I gratefully acknowledge the CEMRACS 2007 organisators who permitted me to stay in the CIRM in Marseilles
and benefit its rich bibliographical ressources.
I also wish to express my gratitude to Yvon Maday and Fr\'ed\'eric Coquel
for fruitful discussions.


\bibliographystyle{plain}

 

 
\end{document}